\documentclass[11pt]{article}
\usepackage[dvips]{graphicx}
\usepackage{amssymb,amsmath,amsthm}
\begin{document}

\vspace{5cm}
\begin{center}
{\Large {\bf BIFURCATION AND STABILITY ANALYSIS OF BISTABLE NEUROMODULES}} \\
\end{center} 

\normalsize
\begin{center}
Stephen Lynch$^*$ and Jon Borresen \\
{\it School of Computing, Mathematics and Digital Technology \\
Manchester Metropolitan University, \\
Manchester, M1 5GD, UK \\
Email: s.lynch@mmu.ac.uk}\\
\end{center}

\noindent
{\bf Abstract}\\

\noindent
This paper presents a stability analysis of simple neuromodules displaying fold bifurcations (leading to hysteresis), 
flip bifurcations (period doubling and undoubling to and from chaos) and Neimark-Sacker bifurcations (quasiperiodic 
and periodic bifurcations). For the first time, bifurcation diagrams are plotted using a feedback mechanism. It is 
shown that the stability curves and bifurcation diagrams must be dealt with simultaneously in order to fully 
understand the dynamics of the systems involved. Synaptic weights, biases and gradients of transfer functions are 
varied and the system is shown to be history dependent. The work can be applied to artificial neural networks and developing 
brains and gives a very important generalization of previous work in this field.\\

\\

$^*$ Corresponding author.\\

E-mail: s.lynch@mmu.ac.uk

\normalsize
\vspace{2cm}
\noindent
{\bf 1. Introduction}

\vspace{0.2cm}
\noindent
In recent times, the disciplines of neural networks and dynamical systems have increasingly coalesced and a new branch %%@
of science called neurodynamics is emerging. This paper is primarily concerned with a simple nonlinear two-neuron %%@
module displaying bistable behaviour, but other nonlinear phenomena have to be discussed in order to fully understand %%@
the dynamics involved. The paper deals with discrete models, but the dynamics of continuous and stochastic models are %%@
important and will now be briefly discussed. 

The basic model of a cell membrane is that of a resistor and capacitor in parallel. The differential equations used to %%@
model the membrane are a variation of the van der Pol equations and represent a simple continuous system. The famous %%@
Fitzhugh-Nagumo oscillator \cite{Fitz, Nagu}, used to model the action potential of a neuron, is a two-variable simplification of %%@
the Hodgkin-Huxley equations \cite{Hodg}. The integrate and fire neuron modelling the slow collection and fast release of %%@
voltage is a biologically inspired model of a neuron and matches the qualitative behaviour very well.
The authors are currently investigating new computer architectures using such coupled oscillator dynamics \cite{Bor1}--\cite{Ash3}. This %%@
recent work on coupled Hodgkin-Huxley equations with time delayed synaptic coupling, demonstrates controlled switching %%@
between synchronized cluster states, leading to reliable logic based programming derived from a biological model. This %%@
work should have significant implications for understanding neural encoding with applications in novel computation. In %%@
1982, Hopfield \cite{Hopf} showed how an analog electrical circuit could behave as a small network of neurons with graded %%@
response. He derived a suitable Lyapunov function for the network to check for stability, and used it as a %%@
content-addressable memory. 

More recently, the authors have had UK and international patents published on binary oscillator computing \cite{Lyn1, Lyn2}.
Arithmetic logic and memory have been demonstarted using Fitzhugh-Nagumo threshold oscillators to perform
binary logic \cite{Bor2}. They are currently seeking industrial partners to investigate possible technological implementations
using Josephson junctions as the threshold oscillators.

It is well known that time continuous systems are often investigated in terms of their discrete time analogues. Using %%@
approximations and algorithms from numerical analysis, a discretized version of the differential equations used to %%@
model neuron dynamics can be computed \cite{Hirs}. One popular method for producing discretized versions of higher dimensional %%@
continuous systems involves constructing Poincar\'{e} maps, see \cite{Lyn3}--\cite{Lyn5}, for example. 

The mechanism of stochastic resonance was introduced in the early 1980's by Benzi et al. \cite{Benz}. The nonlinear 
effect of stochastic resonance at first appears to be counterintuitive---a random noise can make a system 
sensitive to an otherwise undetectable signal. In practice, there is normally an intermediate noise level which 
makes the signal transmission optimal \cite{Fauv}. Under certain conditions, noise can improve signal transmission 
properties of neuronal systems \cite{Tuck}. In general, a two-state model is enough to exhibit the main features of 
stochastic resonance. Zaikin et al.~\cite{Zaik} show that noise induces jumps between different steady states in a 
bistable system. Drover and Ermentrout \cite{Drov} analyzed a nonlinearly coupled system of bistable differential 
equations and showed that there are jumps between critical points and limit cycles. Lynch and Christopher \cite{Lyn6}
have demonstrated multiple bistable loops for certain highly nonlinear differential equations used to model 
wing rock in modern aircraft and surge in jet engines. For these systems, there are jumps from limit cycles to 
limit cycles as certain parameters vary. 

Chaotic and periodic behaviour has been observed in the human brain, see \cite{Duke}--\cite{Pant}, for example. In terms of %%@
mathematical models, Pasemann and his research group have published a number of papers on simple neuromodules %%@
analysing their bistability, chaos, chaos control, synchronization, periodicity and quasiperiodicity, \cite{Pas1}--\cite{Pas5}.
Yuan et al. \cite{Yuan}, study the stability and bifurcation analysis of a discrete-time network for a single fixed point in %%@
terms of Neimark-Sacker bifurcations with special classes of transfer function. They list parameter values where the %%@
fixed point is stable and using the normal form method and the centre manifold theory for discrete systems developed %%@
by Kuznetzov \cite{Kuzn}, they give an algorithm for determining whether a bifurcation from the fixed point is subcritical or %%@
supercritical. In 2008, Guo et al. \cite{Guo} extended this work and considered a simple discrete two-neuron network with %%@
three delays. They derive necessary and sufficient conditions for the asymptotic stability of a fixed point and %%@
investigate the existence of three types of bifurcation, namely, flip bifurcation, fold bifurcation and Neimark-Sacker %%@
bifurcation. Redmond et al. \cite{Redm}, consider a general class of first-order nonlinear delay-differential equations with %%@
reflectional symmetry and apply the theory to simple models from neural networks and atmospheric physics. Their %%@
analysis reveals a Hopf bifurcation curve terminating on a pitchfork bifurcation line at a codimension-two %%@
Takens--Bogdanov point in parameter space. X. Xu \cite{Xu}, investigates local and global Hopf bifurcation in a two-neuron %%@
module with two distinct delays. In all of the papers above, the analysis is restricted to a single fixed point.  

It has long been known that simple neuronal models can display steady-state behaviour, periodic behaviour, %%@
quasiperiodic behaviour and chaos. One major factor that has been overlooked in studying simple neuromodules is %%@
feedback. The two essential ingredients for bistability, or hysteresis, are nonlinearity and feedback \cite{Lyn1}--\cite{Lyn3}. Both %%@
of these processes are inherently present in neuronal dynamics. Bistable behaviour often occurs in physical systems %%@
for only a narrow range of parameter values and hence can be difficult to observe. Mathematical analysis can be used %%@
to test for stability, and bistable and unstable regions can be plotted for a range of parameter values. Bistable %%@
phenomena have been proven to exist in a wide range of physical forms. They are present in biological systems \cite{Gaver,Lyn7}, %%@
human brain dynamics \cite{Meye}, economics \cite{Goc}, mechanical systems \cite{Go}, electric circuits \cite{Bor2}, chemical kinetics \cite{Scot}, %%@
astrochemical cloud models \cite{Nejad}, psychology \cite{Fish,Arro}, neural networks \cite{Izhi} and nonlinear optics \cite{Stee,Lyn8}, where optical bistability has potential applications in high speed all-optical signal processing and all-optical computing. Parameters such as %%@
weights, biases and gradients of transfer functions can be altered when an artificial neural network is learning. Biologically, these same variables obviously go through similar changes as the biological brain develops.

Section 2 is concerned with the mathematical model of a single neuron considered by Pasemann \cite{Pas2}, and a simple %%@
two-neuron module is investigated in Section 3. The analysis can be seen to match very well with the numerical results %%@
displayed in bifurcation diagrams incorporating a feedback mechanism, however, there are exceptional cases which have %%@
not been discussed elsewhere. Conclusions and further work will be discussed in Section 4.\\

\noindent
{\bf 2. Stability analysis of a single neuron}
\vspace{0.2cm}

The simple chaotic neuron investigated by Pasemann \cite{Pas2} is modelled using the discrete equation
\begin{equation}
x_{n+1}=b+\gamma x_n+w\sigma(x_n), \label{eq:1}
\end{equation}
where $x_n$ is the activation level at time $n$, $b$ is a bias, $w$ is a self-weight and, $0<\gamma<1$, is a decay %%@
rate depicting the dissipative property of a real neuron. The output is given by the unipolar transfer function, %%@
$\sigma(x)=1/(1+e^{-x})$. This system can display hysteresis and period-doubling bifurcations to and from chaos.
This system cannot display quasiperiodic behaviour. Fixed points of period one satisfy the equation
$$x=\gamma x+b+w \sigma(x),$$
where $x$ is the fixed point of period one. The stability edges for $x$ occur when $|f'(x)|<1$, 
$f(x)=b+\gamma x+w\sigma(x)$. Thus, the fixed point is stable as long as
$$|\gamma+w\sigma'(x)|<1.$$
The parametric equations defining the stability boundaries are given by 
$$B^{+1}: \quad b=(1-\gamma)x-\frac{1-\gamma}{1-\sigma(x)}, \quad w=\frac{1-\gamma}{\sigma'(x)},$$
$$B^{-1}: \quad b=(1-\gamma)x+\frac{1+\gamma}{1-\sigma(x)}, \quad w=-\frac{1+\gamma}{\sigma'(x)},$$
where $B^{+1}$ defines the bistable boundary, where fixed points undergo fold bifurcations. The curve depicted by %%@
$B^{-1}$ determines where the system goes unstable, and fixed points undergo flip bifurcations. The 
parametric curves are plotted in Figure 1 when $\gamma=0.5$. In this case, the bistable region is isolated from 
instabilities.

\begin{figure}[h!]
\centerline{\includegraphics[width=80mm]{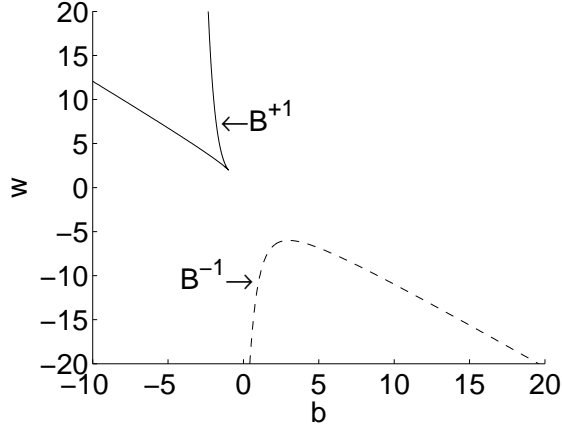}} 
\caption{\small{Stability diagram for the nonlinear neuron when $\gamma=0.5$. The curve $B^{+1}$ is the 
boundary for bistable operation and the curve $B^{-1}$ is the border for period-doubling bifurcations.}}
\end{figure} 

Bistable behaviour for the single neuron is demonstrated in Figures 2(a) and 2(b). In Figure 2(a), the bias $b$ is %%@
ramped up from $b=-5$ to $b=5$, and then ramped down again. One parameter is varied and the solution to the previous %%@
iterate is used as the initial condition for the next iterate; in this way, a feedback mechanism is introduced. In %%@
this case, there is a history associated with the process and only one point is plotted for each value of the 
parameter. This method has been labelled as the second iterative method by Lynch \cite{Lyn3}--\cite{Lyn5}. The bistable region 
matches with the results of the linear stability analysis given in Figure 1. In Figure 2(b), a larger bistable region %%@
is produced as both parameters, $b$ and $w$, are increased and then decreased simultaneously. It should be pointed out %%@
at this stage that in most physical situations many parameters would be varying simultaneously and it is only in rare %%@
cases that most parameters would be kept constant as one parameter is varied.  

\begin{figure}[htbp]
\centerline{(a)\includegraphics[width=75mm]{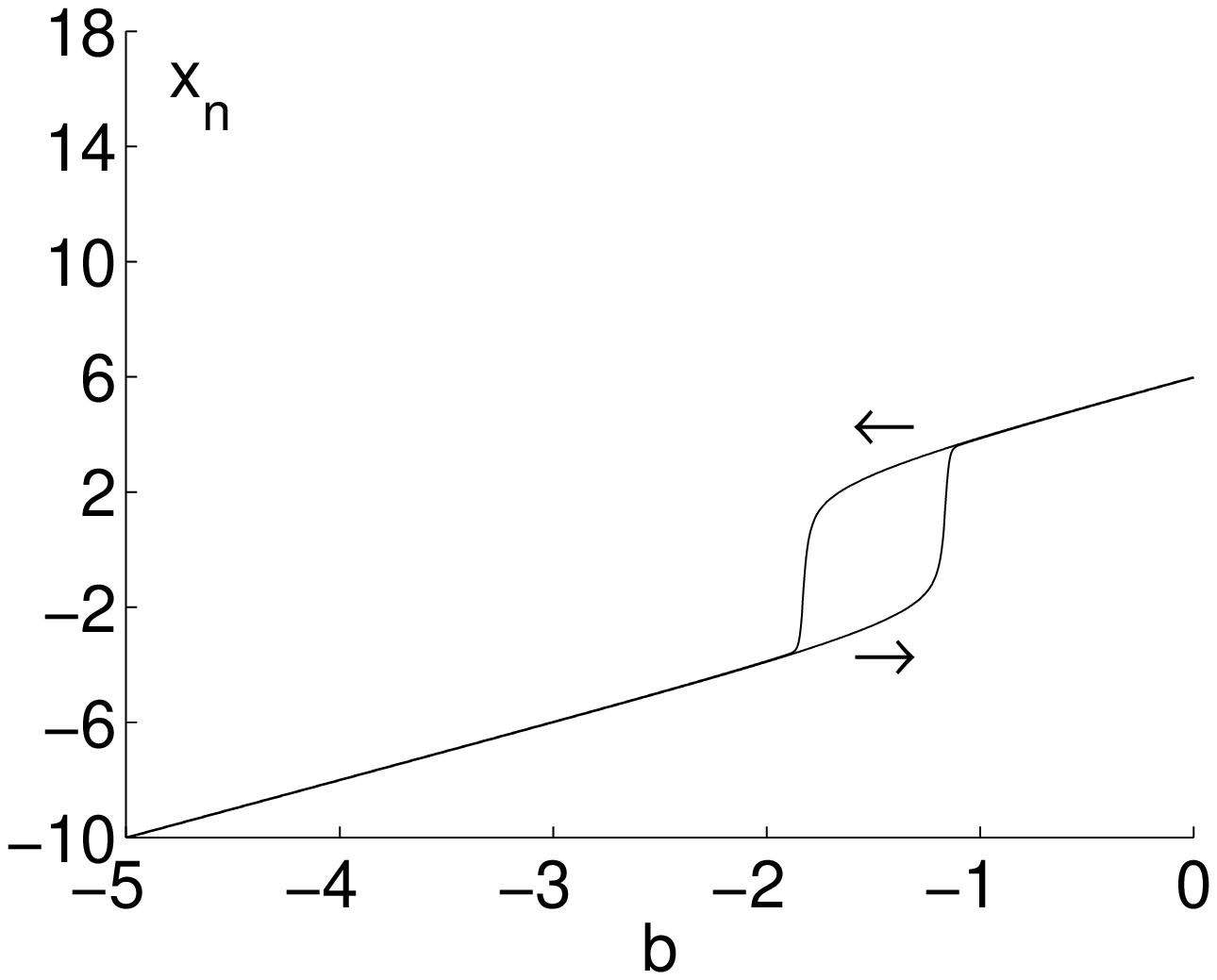} 
(b)\includegraphics[width=75mm]{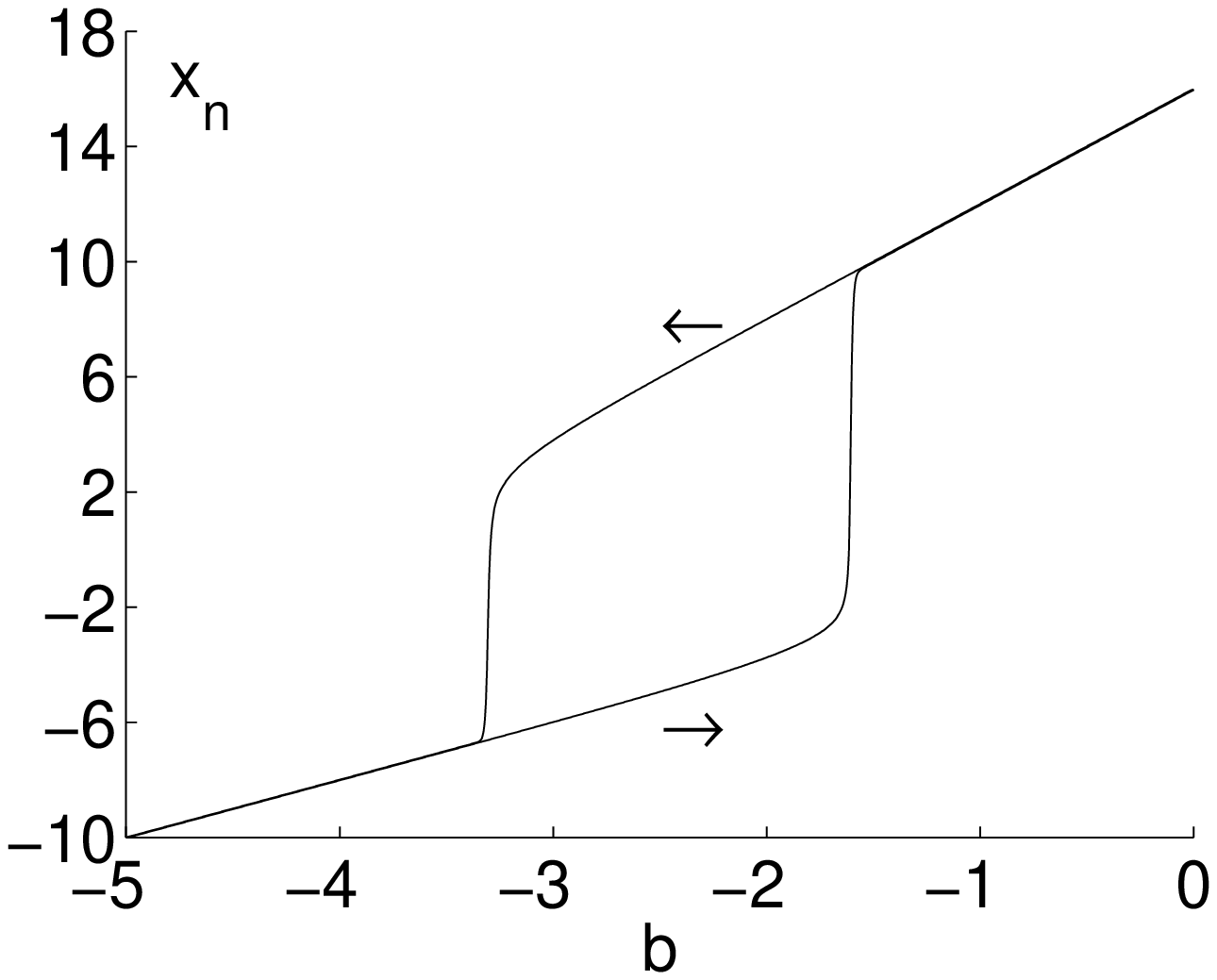}} 
\caption{\small{Bifurcation diagrams for system (\ref{eq:1}) when $\gamma=0.5$ using a stepsize of $0.01$. (a) The %%@
bias $b$ is ramped up from $-5$ to $0$ and then ramped down again, the parameter $w$ is fixed at $w=3$. (b) The bias %%@
$b$ and weight $w$ are ramped up and down simultaneously, $-5\leq b\leq 0$ and $3\leq w\leq 5$.}}
\end{figure}

Figures 3(a) and 3(b) show bifurcation diagrams for system (\ref{eq:1}), when $\gamma=0.5$ and the parameter $b$ is %%@
ramped up and down, $(-5\leq b\leq 5)$, and parameter $w$ is simultaneously ramped down then up, $(-10\leq w \leq %%@
10)$. In Figure 3(a), the initial condition chosen at $\gamma=0.5, b=-5, w=10$ was $x(0)=-10$, however, in Figure 3(b) %%@
the initial condition chosen for the same set of parameter values was $x(0)=10$. In Figure 3(a), there is a small %%@
bistable region for $b$ close to $-1$, which shrinks as the step size used in the iterative scheme is decreased. For %%@
increasing values of the parameter $b$, the solution remains on the steady-state even though the fixed point becomes %%@
unstable at $b \approx 3$. As the parameter $b$ is decreased, from $b=5$, the solution becomes unstable (a stable %%@
period two orbit) and then returns to the steady-state at $b \approx 3$, once more, where it remains as $b$ is ramped %%@
back down to $b=-5$ and $w$ is ramped back up to $w=10$. In Figure 3(b), there is a large bistable region for %%@
$-5<b<-1$, approximately. For increasing values of the parameter $b$, the solution again remains on the steady-state %%@
on the upper branch of the hysteresis cycle. As the parameter $b$ is decreased, from $b=5$, the solution becomes %%@
unstable (period two) and then returns to the steady-state at $b \approx 3$, however, as the parameter $b$ is %%@
decreased further, the steady-state now follows the lower branch of the hysteresis cycle, corresponding to the other %%@
stable fixed point of the system. Hence it has been shown that the bifurcation diagrams are sensitive to both step %%@
size and the choice of initial conditions. It is also evident that it is possible for the steady-state to remain on a %%@
path even though that path is unstable, as Figures 3(a) and 3(b) demonstrate.\\

\begin{figure}[h!]
\centerline{\includegraphics[width=75mm]{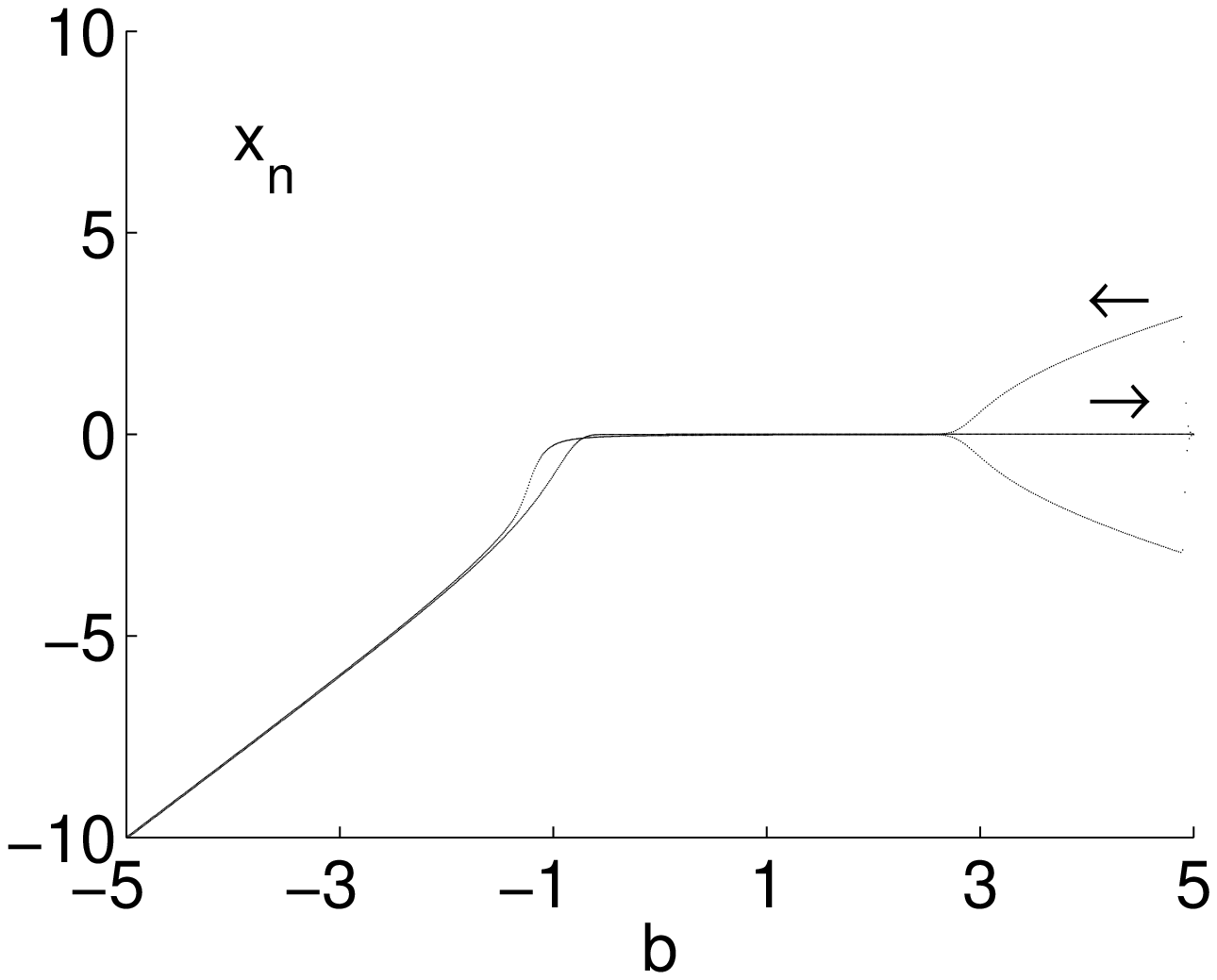} \\
\includegraphics[width=75mm]{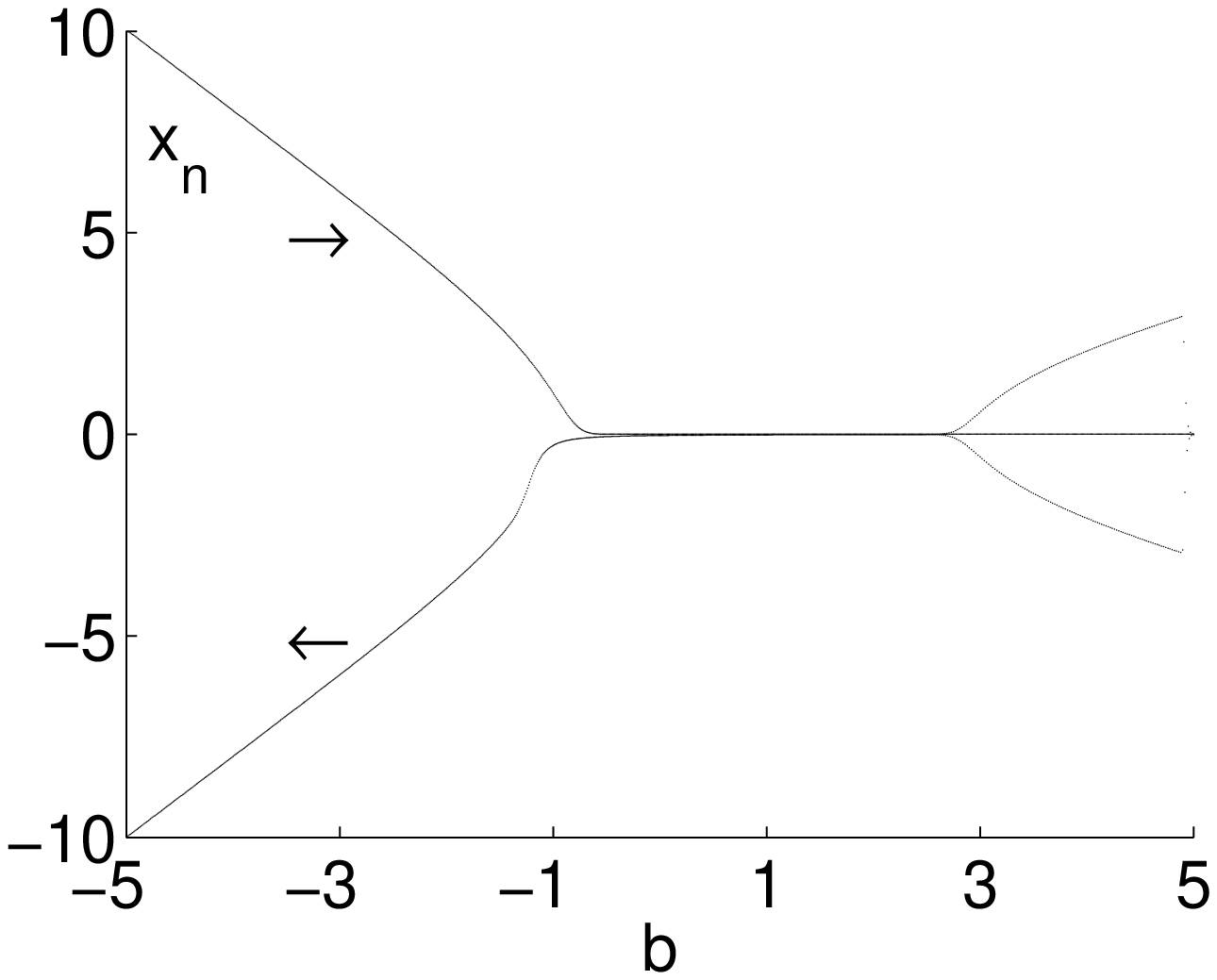}} 
\begin{tabbing}
\hspace*{3cm} \= \hspace*{8cm} \= \hspace*{2cm} \kill
 \> (a) \> (b) 
\end{tabbing}
\caption{\small{Bifurcation diagrams for system (\ref{eq:1}) when $\gamma=0.5$ using a stepsize of $0.01$. Parameters %%@
$b$ and $w$ are varied simultaneously, $b$ is ramped up then down and parameter $w$ is ramped down then up, $(-5\leq %%@
b\leq 5,-10\leq w \leq 10)$. (a) The initial condition chosen at $b=-5, w=10$ is $x(0)=-10$. (b) The initial condition %%@
for the same set of parameter values is $x(0)=10$.}}
\end{figure}

The next section is concerned with a simple chaotic two-neuron module that can display hysteresis and periodic and %%@
quasiperiodic behaviour. \\

\noindent
{\bf 3. Stability analysis of a simple two-neuron module}
\vspace{0.2cm}

Consider the simple two-neuron module as depicted in Figure 4. The discrete dynamical system used to model the %%@
neuromodule is given by
\begin{equation}
x_{n+1}=b_1+w_{11}\sigma_1\left(x_n\right)+w_{12}\sigma_2\left(y_n\right), \quad
y_{n+1}=b_2+w_{21}\sigma_1\left(x_n\right)+w_{22}\sigma_2\left(y_n\right), \label{eq:2}
\end{equation}
where $x_n$, $y_n$ are the activation levels of the neurons; $b_1$, $b_2$ are biases; $w_{ij}$ are synaptic weights %%@
and $\sigma_1(x)$, $\sigma_2(y)$ are transfer functions. This simple system displays chaos, periodicity and %%@
quasiperiodicity, characteristics that are evident in much more complicated neuromodules. Without loss of generality, %%@
assume in the work to follow that the weight $w_{22}$ is zero and the transfer functions are defined by the bipolar %%@
functions $\sigma_1(x)=\tanh(\alpha x)$, $\sigma_2(y)=\tanh(\beta y)$, this will make the analysis more manageable.

\begin{figure}[h!]
\centerline{\includegraphics[width=30mm]{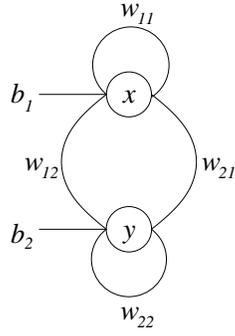}} 
\caption{\small{A recurrent two-neuron module.}}
\end{figure}

Fixed points of period one, which will correspond to steady-state behaviour, are found by solving the equations,
$x_{n+1}=x_n=x$, say, and $y_{n+1}=y_n=y$, say. Hence from system (\ref{eq:2}),
\begin{equation}
b_1=x-w_{11}\tanh(\alpha x)-w_{12}\tanh(\beta y),  \label{eq:3}
\end{equation}
and
\begin{equation}
y=b_2+w_{21}\tanh(\alpha x).  \label{eq:4}
\end{equation}
The stability of these fixed points is determined by a linearization technique [13--15]. Each fixed point can be transformed to the origin and the nonlinear terms can be discarded after taking a Taylor series expansion. The Jacobian matrix is given by
$$J=\left( \begin{array}{cc}
					\frac{\partial P}{\partial x} & \frac{\partial P}{\partial y} \\
					\frac{\partial Q}{\partial x} &\frac{\partial Q}{\partial y}
					\end{array}
                   \right),$$
where $P\left(x_n,y_n\right)=x_{n+1}$,  $Q\left(x_n,y_n\right)=y_{n+1}$. Therefore
$$J=\left( \begin{array}{cc}
					\alpha w_{11}\mathrm{sech}^2(\alpha x) & \beta w_{12}\mathrm{sech}^2(\beta y) \\
					\alpha w_{12}\mathrm{sech}^2(\alpha x) & 0
					\end{array}
                   \right).$$				  
The stability conditions for the fixed points are determined by considering the eigenvalues, trace and determinant of %%@
the Jacobian matrix. The characteristic equation is given by
\begin{equation}
\lambda^2-\alpha w_{11} \mathrm{sech}^2(\alpha x)\lambda-\alpha \beta w_{12}w_{21}\mathrm{sech}^2(\alpha %%@
x)\mathrm{sech}^2(\beta y)=0. \label{eq:5}
\end{equation}
The fixed points undergo a fold bifurcation when $\lambda=+1$. In this case, equation (\ref{eq:5}) gives
\begin{equation}
w_{12}=\frac{1-\alpha w_{11}\mathrm{sech}^2(\alpha x)}{\alpha \beta w_{21}\mathrm{sech}^2(\alpha %%@
x)\mathrm{sech}^2(\beta y)}. \label{eq:6} 				   
\end{equation}
The fixed points undergo a flip bifurcation when $\lambda=-1$. In this case, equation (\ref{eq:5}) gives
\begin{equation}
w_{12}=\frac{1+\alpha w_{11}\mathrm{sech}^2(\alpha x)}{\alpha \beta w_{21}\mathrm{sech}^2(\alpha %%@
x)\mathrm{sech}^2(\beta y)}. \label{eq:7} 				   
\end{equation}
The fixed points undergo a Neimark-Sacker bifurcation when $\det(J)=1$ and $|tr(j)|<2$. Thus
\begin{equation}
w_{12}=-\frac{1}{\alpha \beta w_{21}\mathrm{sech}^2(\alpha x)\mathrm{sech}^2(\beta y)}. \label{eq:8} 				   
\end{equation}
\begin{figure}[h!]
\centerline{\includegraphics[width=60mm]{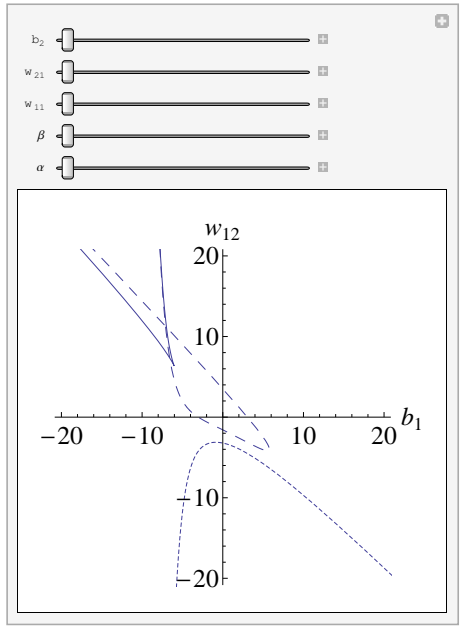}} 
\caption{\small{Mathematica parameter slider generated with the Manipulate command. The user can change parameters, %%@
$\alpha$, $\beta$, $b_2$, $w_{11}$ and $w_{21}$ interactively. The solid curve represents the boundary for bistable %%@
operation and will be labelled as $B^{+1}$ from now on; the dashed curve is the border for period-doubling %%@
bifurcations and will be labelled $B^{-1}$ from now on, and the dotted curve represents the boundary for %%@
Neimark-Sacker bifurcations, and will be labelled $B^{NS}$, from now on.}}
\end{figure}

Figure 5 shows a typical stability diagram for system (\ref{eq:2}) using equations (\ref{eq:3}) and (\ref{eq:4}) and equations (\ref{eq:6})---(\ref{eq:8}), and is %%@
generated using the Manipulate command in Mathematica [14]. As the slider buttons are moved left and right, the stability diagram alters appropriately in an animated sequence.\\

There are three main cases to consider for system (\ref{eq:2}):\\

\begin{figure}[h!]
\centerline{\includegraphics[width=80mm]{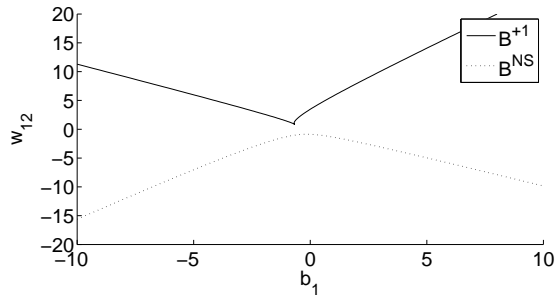}} 
\caption{\small{A typical stability diagram when $w_{11}=0$. In this case, $b_2=3$, $w_{21}=5$, $\alpha=1$ and %%@
$\beta=0.3$.}}
\end{figure}

\noindent
CASE 1: Suppose that the parameter $w_{11}=0$, so neuron $x$ has no self-connection. \\

A typical stability diagram is plotted in Figure 6. It can be seen that the bistable curve $B^{+1}$ is isolated from %%@
the Neimark-Sacker stability curve $B^{NS}$, say. Two possible bifurcation diagrams are plotted in Figure 7. In Figure %%@
7(a), the parameter $b_1$ is ramped up from $b_1=-5$ up to $b_1=5$ and then ramped down again for $w_{12}=5$ and the %%@
parameter values listed in the figure. It can be seen that a large counterclockwise hysteresis loop is present for %%@
$-4<b_1<1$, approximately, which is in agreement with the bounding bistability curves plotted in Figure 6.  In Figure %%@
7(b), the parameter $b_1$ is ramped up from $b_1=-5$ up to $b_1=5$ and then ramped down again for $w_{12}=-4$ and the %%@
parameter values listed in the figure. It can be seen that a large period-2 loop is present for $-3<b_1<4$, %%@
approximately, which again agrees with the bounding Neimark-Sacker curve plotted in Figure 6. \\

\begin{figure}[h!]
\centerline{\includegraphics[width=75mm]{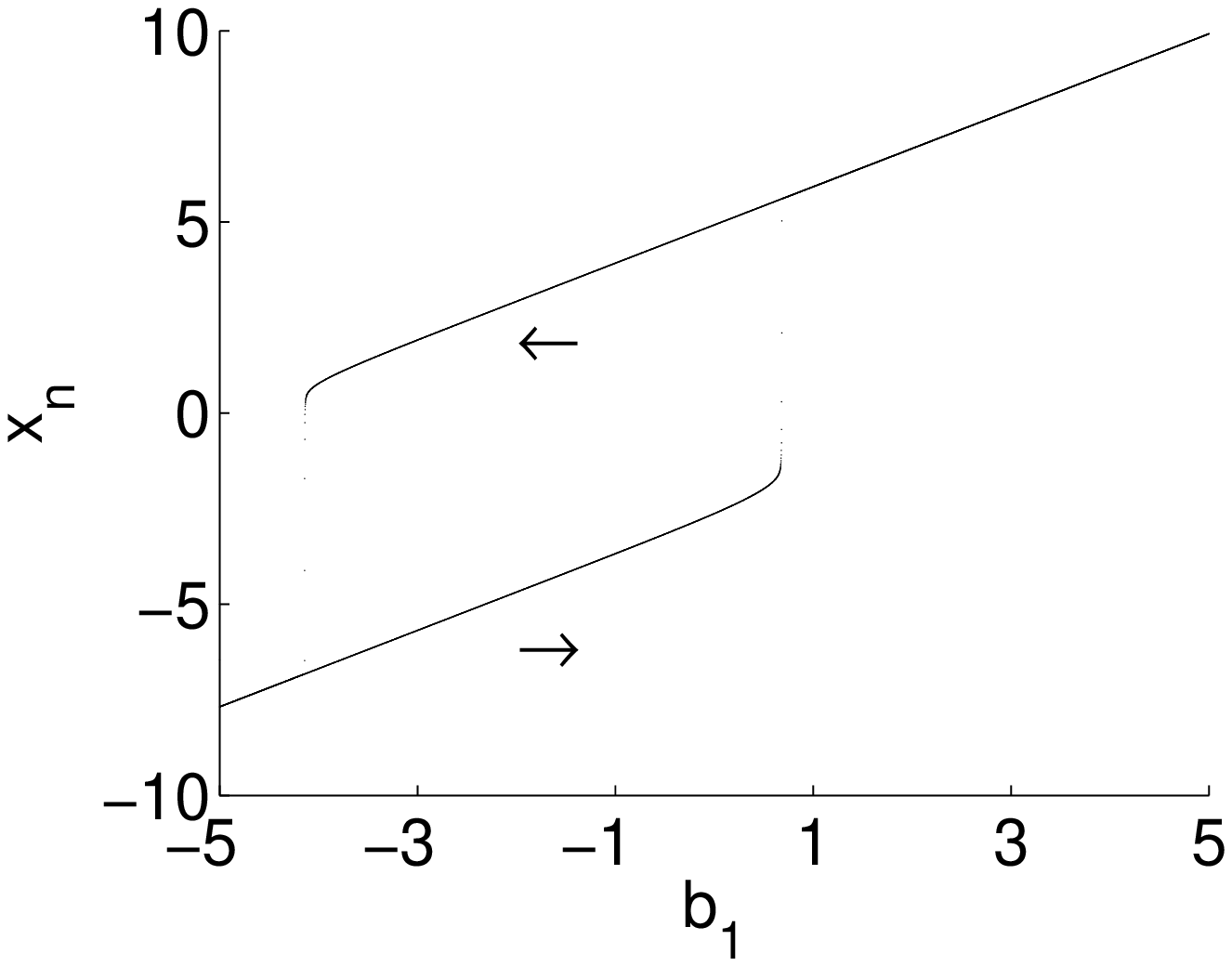} \hspace{0.2cm} 
\includegraphics[width=75mm]{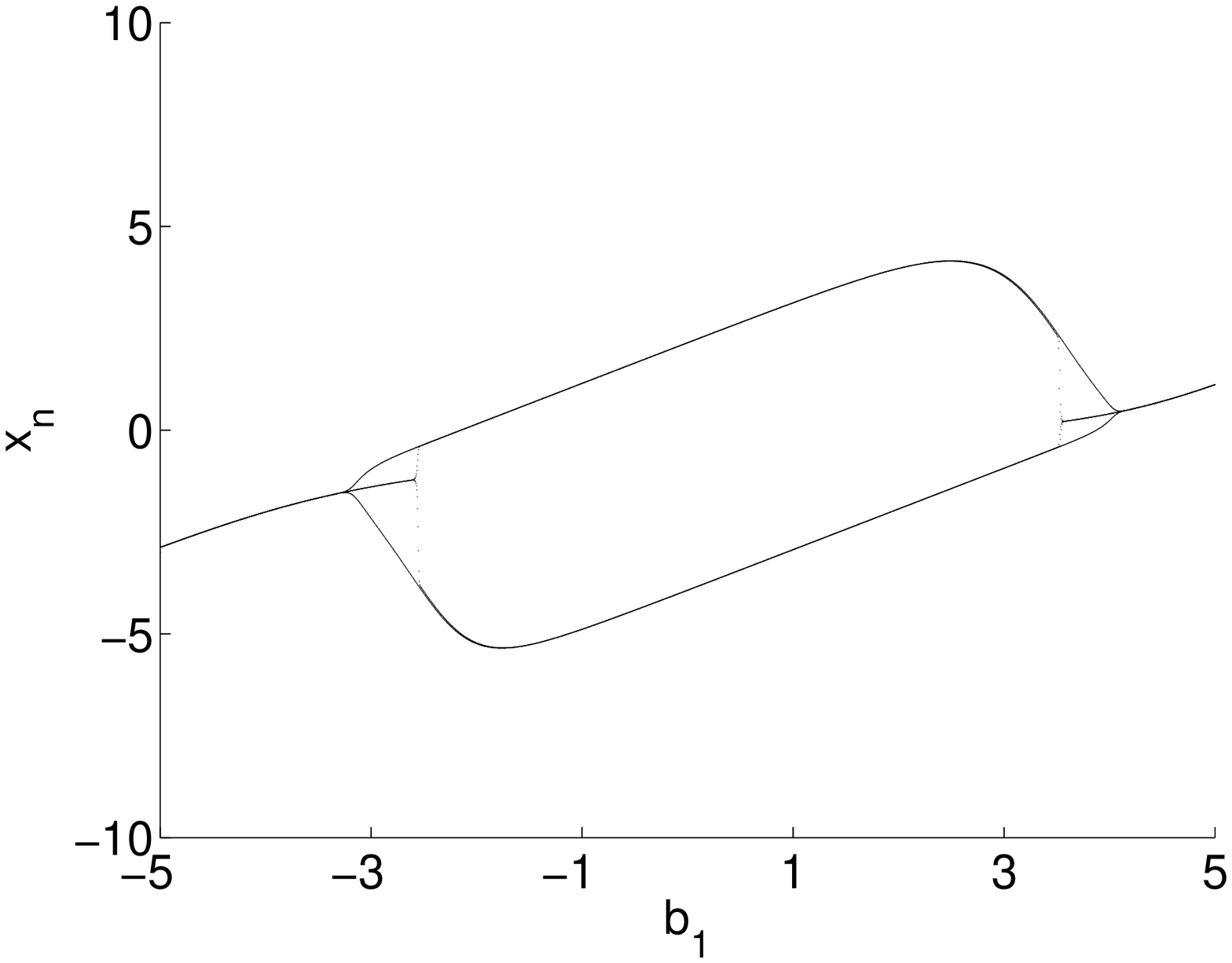}} 
\begin{tabbing}
\hspace*{3cm} \= \hspace*{8cm} \= \hspace*{2cm} \kill
 \> (a) \> (b) 
\end{tabbing}
\caption{\small{Typical bifurcation diagrams for system (\ref{eq:2}) when $w_{11}=0$, $b_2=3$, $w_{21}=5$, $\alpha=1$ %%@
and $\beta=0.3$ using a stepsize of $0.01$. The parameter $b_1$ is ramped up then down, $-5\leq b_1\leq 5$. (a) The %%@
parameter $w_{12}=5$ and the initial conditions are $x(0)=-7$, $y(0)=-2$. (b) The parameter $w_{12}=-4$ and the %%@
initial conditions are $x(0)=-3$, $y(0)=-2$.}}
\end{figure}

\begin{figure}[h!]
\centerline{\includegraphics[width=80mm]{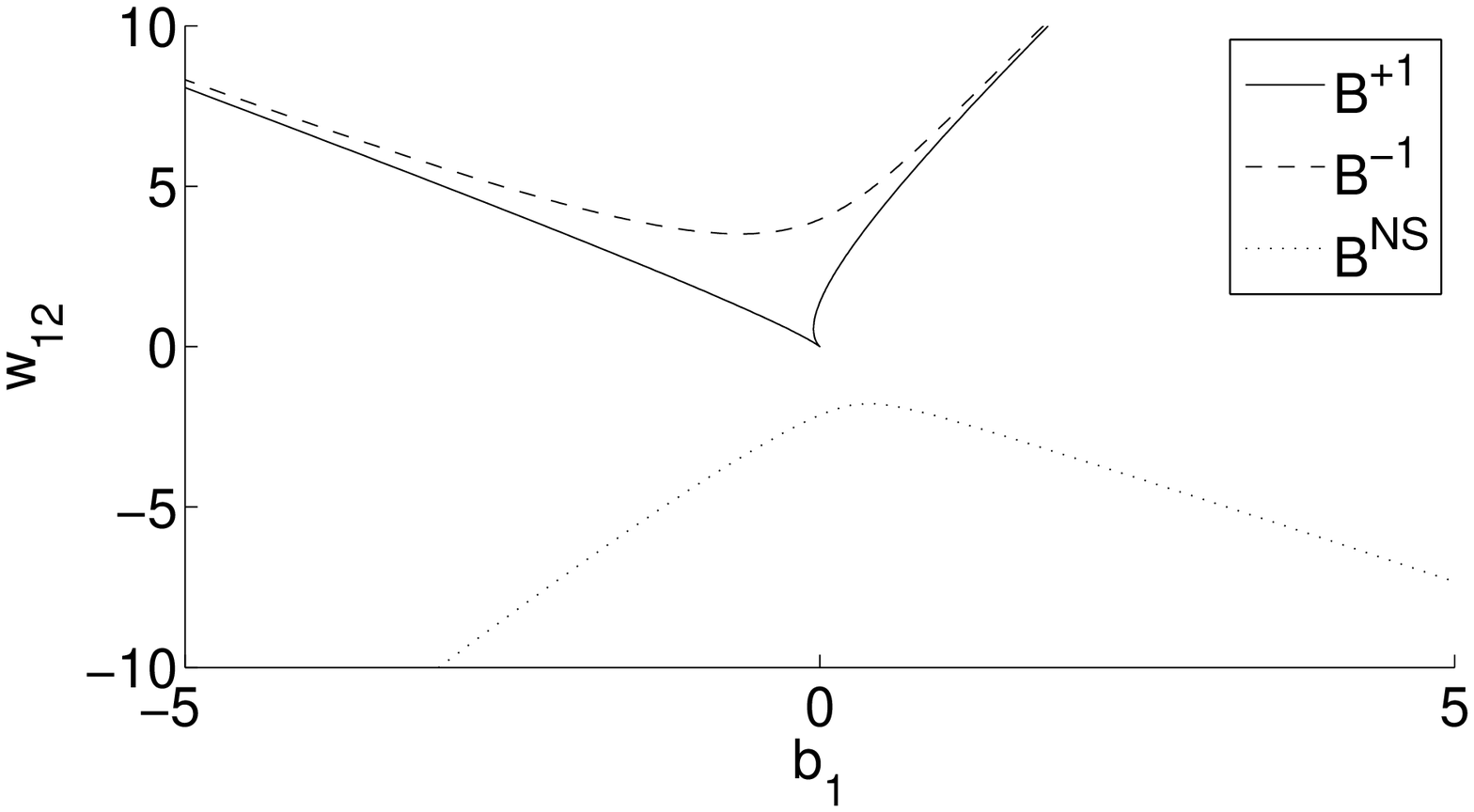}} 
\caption{\small{A typical stability diagram when $w_{11}>0$. In this case, $w_{11}=1$, $b_2=1$, $w_{21}=2$, $\alpha=1$ %%@
and $\beta=0.3$.}}
\end{figure}

\begin{figure}[h!]
\centerline{\includegraphics[width=75mm]{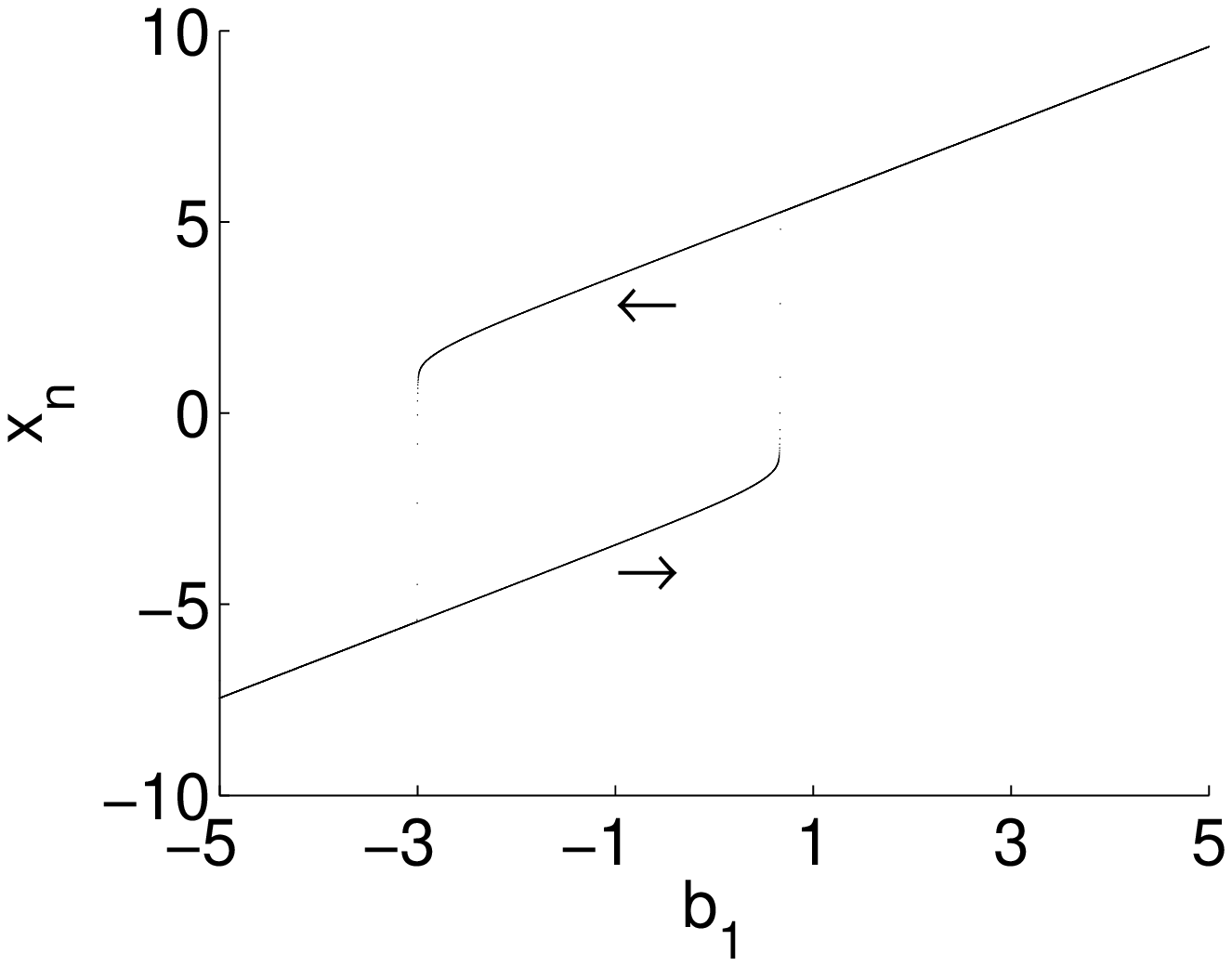} \hspace{0.2cm} 
\includegraphics[width=75mm]{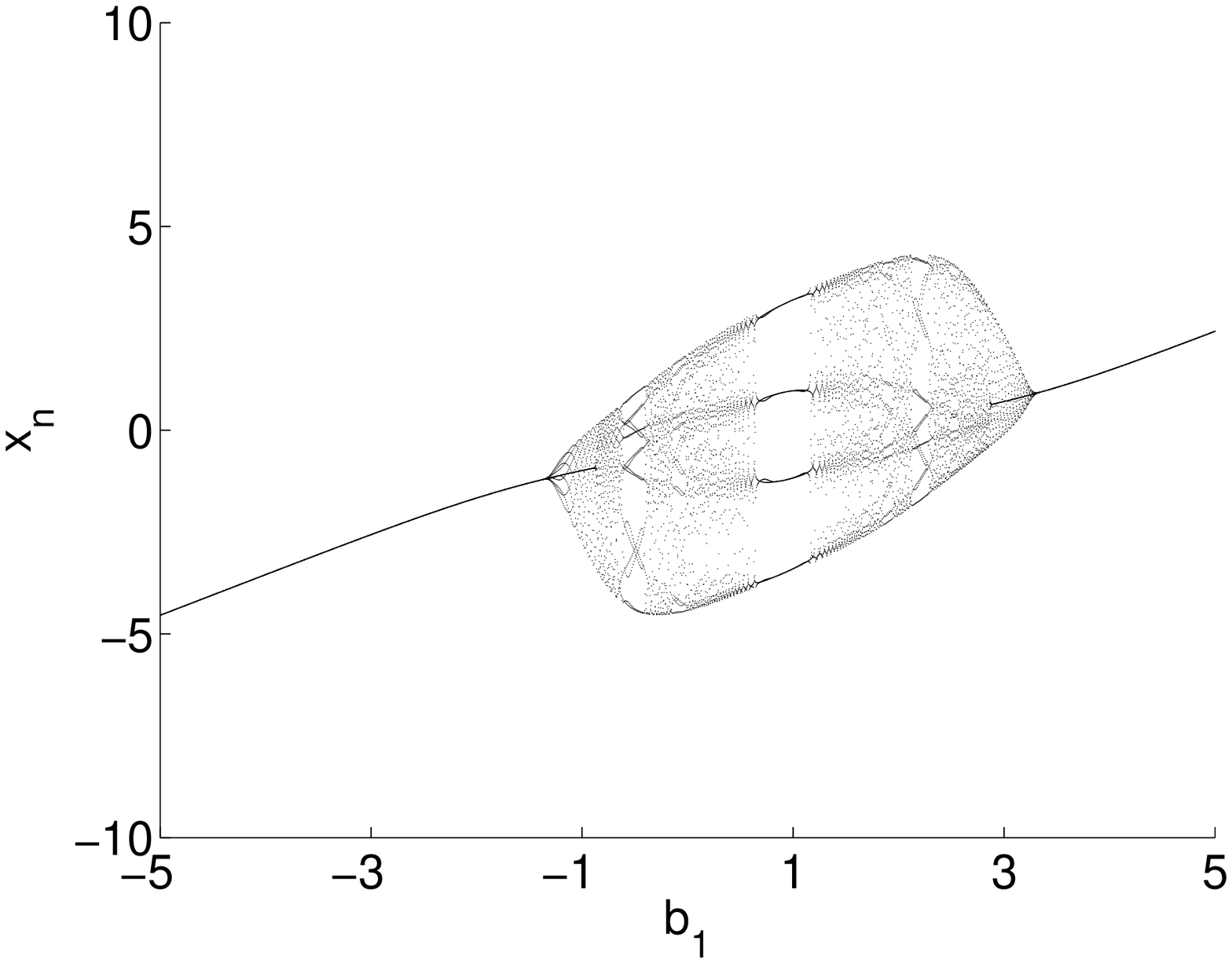}} 
\begin{tabbing}
\hspace*{3cm} \= \hspace*{8cm} \= \hspace*{2cm} \kill
 \> (a) \> (b) 
\end{tabbing}
\caption{\small{Typical bifurcation diagrams for system (\ref{eq:2}) when $w_{11}>0$. In this case, $w_{11}=1$, %%@
$b_2=1$, $w_{21}=2$, $\alpha=1$ and $\beta=0.3$. The parameter $b_1$ is ramped up then down, $-5\leq b_1\leq 5$. (a) %%@
The parameter $w_{12}=5$ and the initial conditions are $x(0)=-7$, $y(0)=-1$. (b) The parameter $w_{12}=-5$ and the %%@
initial conditions are $x(0)=-5$, $y(0)=-1$.}}
\end{figure}

\begin{figure}[htbp]
\centerline{\includegraphics[width=90mm]{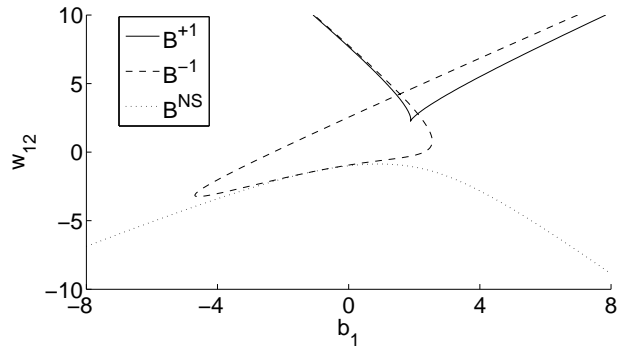}} 
\caption{\small{A typical stability diagram when $w_{11}<0$. In this case, $w_{11}=-2$, $b_2=-3$, $w_{21}=5$, %%@
$\alpha=1$ and $\beta=0.3$.}}
\end{figure}

\noindent
CASE 2: Neuron $x$ is a self-excitatory neuron, where $w_{11}>0$.\\

A typical stability diagram is plotted in Figure 8. The unstable curve bounded by $B^{-1}$ does not affect the %%@
bistable operation when $w_{11}>0$. Two possible bifurcation diagrams are plotted in Figure 9. In Figure 9(a), the %%@
parameter $b_1$ is ramped up from $b_1=-5$ up to $b_1=5$ and then ramped down again for $w_{12}=5$ and the parameter %%@
values listed in the figure. It can be seen that a large counterclockwise hysteresis loop is present for $-3<b_1<1$, %%@
approximately, which is in agreement with the bounding bistability curves plotted in Figure 8.  In Figure 9(b), the %%@
parameter $b_1$ is ramped up from $b_1=-5$ up to $b_1=5$ and then ramped down again for $w_{12}=-5$ and the parameter %%@
values listed in the figure. It can be seen that a large quasiperiodic and periodic loop is present for $-1<b_1<3$, %%@
approximately, which again agrees with the bounding Neimark-Sacker curve plotted in Figure 8.\\

\noindent
CASE 3: Neuron $x$ is a self-inhibatory neuron, where $w_{11}<0$.\\

A typical stability diagram is plotted in Figure 10. The unstable curve bounded by $B^{-1}$ does affect the bistable %%@
operation when $w_{11}<0$. Two possible bifurcation diagrams are plotted in Figure 11. In Figure 11(a), the parameter %%@
$b_1$ is ramped up from $b_1=-8$ up to $b_1=8$ and then ramped down again for $w_{12}=5$ and the parameter values %%@
listed in the figure. It can be seen that a large unstable period-2 cycle has encroached upon the counterclockwise %%@
hysteresis loop for $3<b_1<6$, approximately, which is in agreement with the bounding curves $B^{+1}$ and $B^{-1}$, %%@
plotted in Figure 10.  In Figure 11(b), the parameter $b_1$ is ramped up from $b_1=-5$ up to $b_1=5$ and then ramped %%@
down again for $w_{12}=-5$ and the parameter values listed in the figure. It can be seen that a large period-2 loop is %%@
present for $-3<b_1<0$, approximately, and a quasiperiodic loop is present for $0<b_1<2$, approximately, which again %%@
agrees with the bounding curves $B^{-1}$ and $B^{NS}$, plotted in Figure 10. One can check that the behavior of a discrete system is quasiperiodic by plotting a power spectrum using a fast Fourier transform.

\begin{figure}[htbp]
\centerline{\includegraphics[width=75mm]{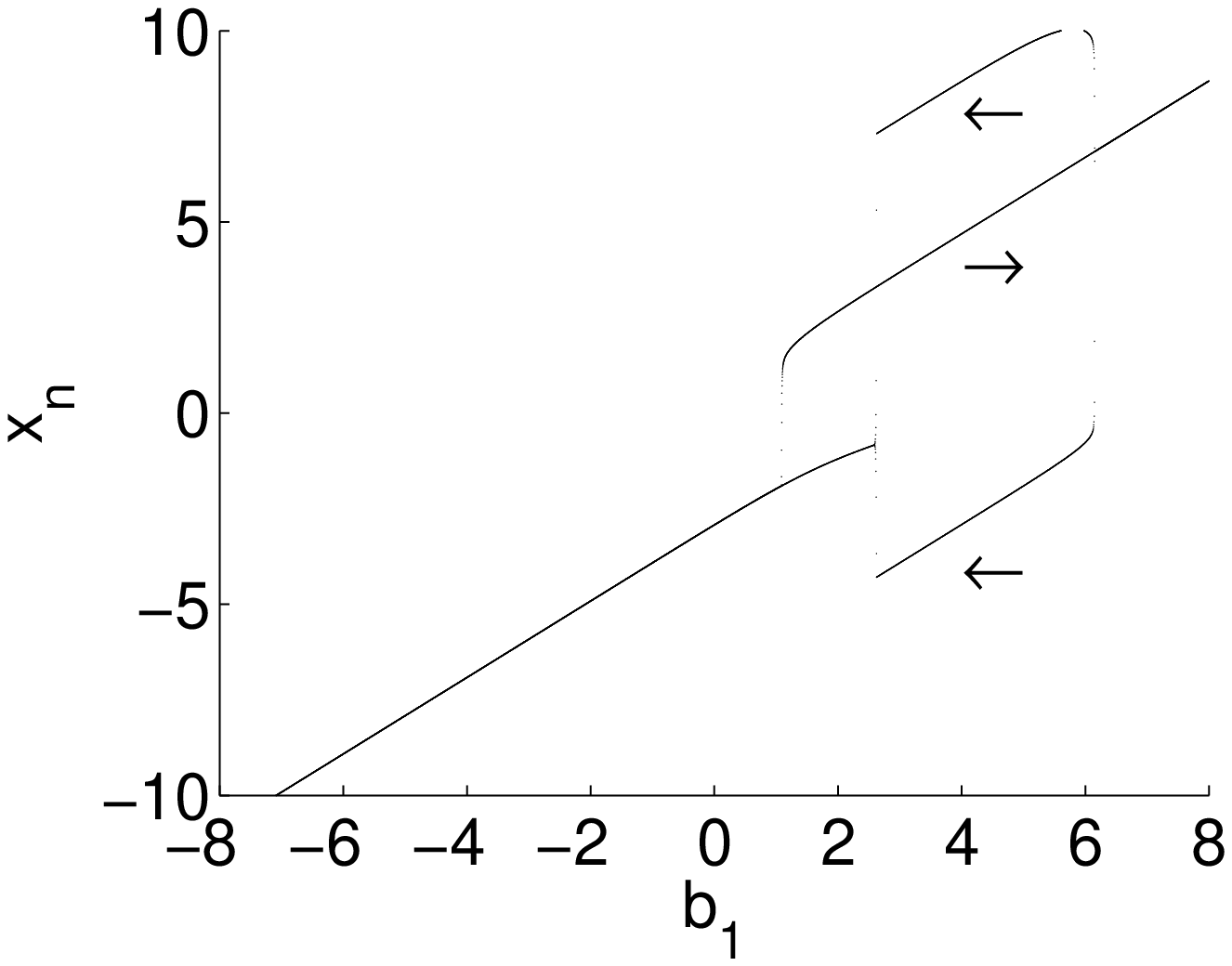} \hspace{0.2cm} 
\includegraphics[width=75mm]{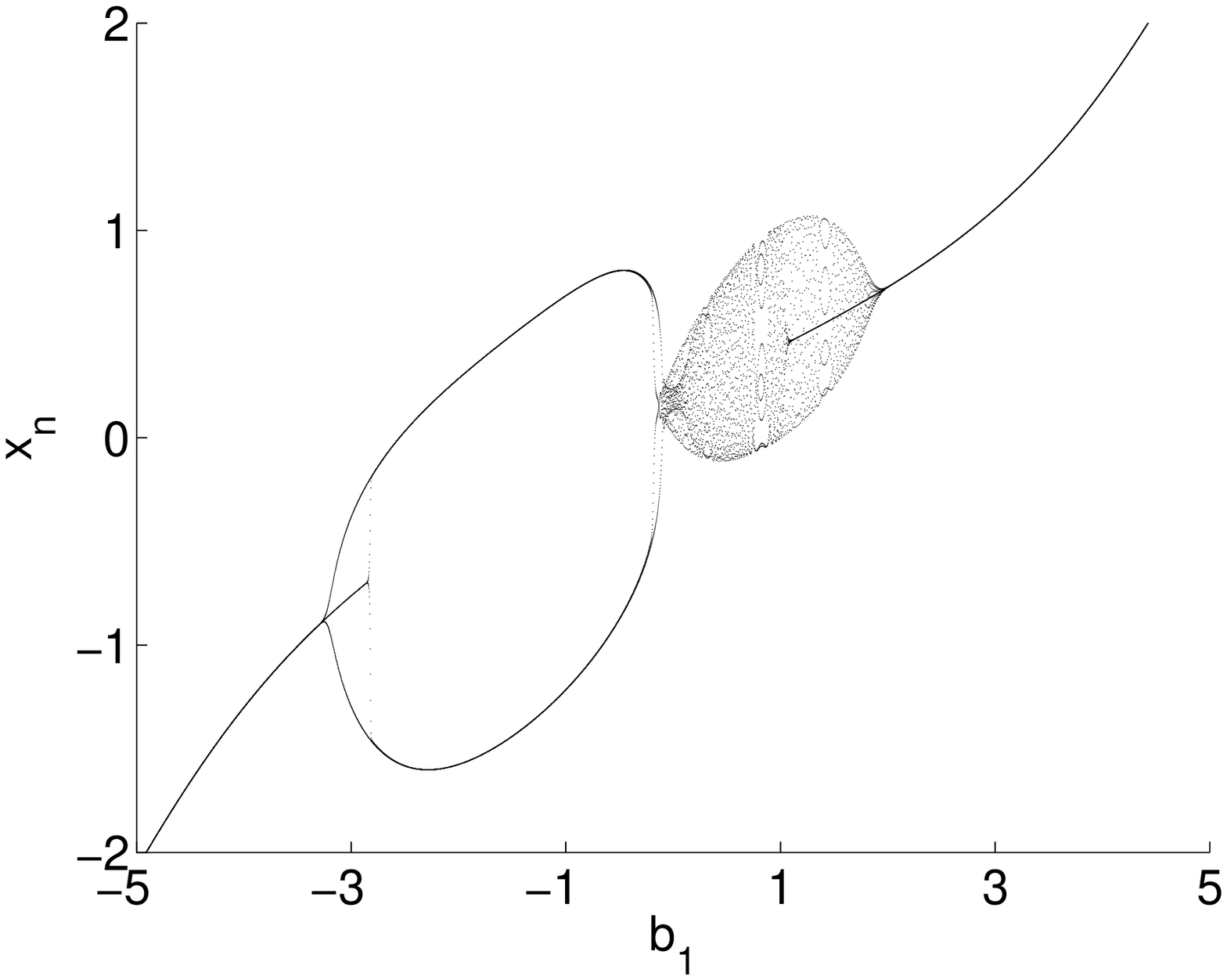}} 
\begin{tabbing}
\hspace*{3cm} \= \hspace*{8cm} \= \hspace*{2cm} \kill
 \> (a) \> (b) 
\end{tabbing}
\caption{\small{Typical bifurcation diagrams for system (\ref{eq:2}) when $w_{11}<0$. In this case, $w_{11}=-2$, %%@
$b_2=-3$, $w_{21}=5$, $\alpha=1$ and $\beta=0.3$. (a) The parameter $b_1$ is ramped up and down, $-8 \leq b_1 \leq 8$, %%@
$w_{12}=5$ and the initial conditions are $x(0)=-11$, $y(0)=-8$. (b)  The parameter $b_1$ is ramped up and down, $-5 %%@
\leq b_1 \leq 5$, $w_{12}=-1$ and the initial conditions are $x(0)=-2$, $y(0)=-8$.}}
\end{figure}

Now fix the parameter $b_1$ and vary the parameter $w_{12}$ for the parameter values listed in Figure 10. The %%@
parameter $b_1$ is fixed at $b_1=1$ and the parameter $w_{12}$ is ramped down then up for $-10 \leq w_{12}\leq 10$. %%@
Figure 12(a) shows the bifurcation diagram obtained when the initial conditions chosen are $x(0)=-7$, $y(0)=-7$ and %%@
Figure 12(b) shows the corresponding bifurcation diagram when the initial conditions are $x(0)=4$ and $y(0)=2$. In %%@
both cases, there is unstable behaviour swamping the bifurcation diagrams.

\begin{figure}[h!]
\centerline{\includegraphics[width=75mm]{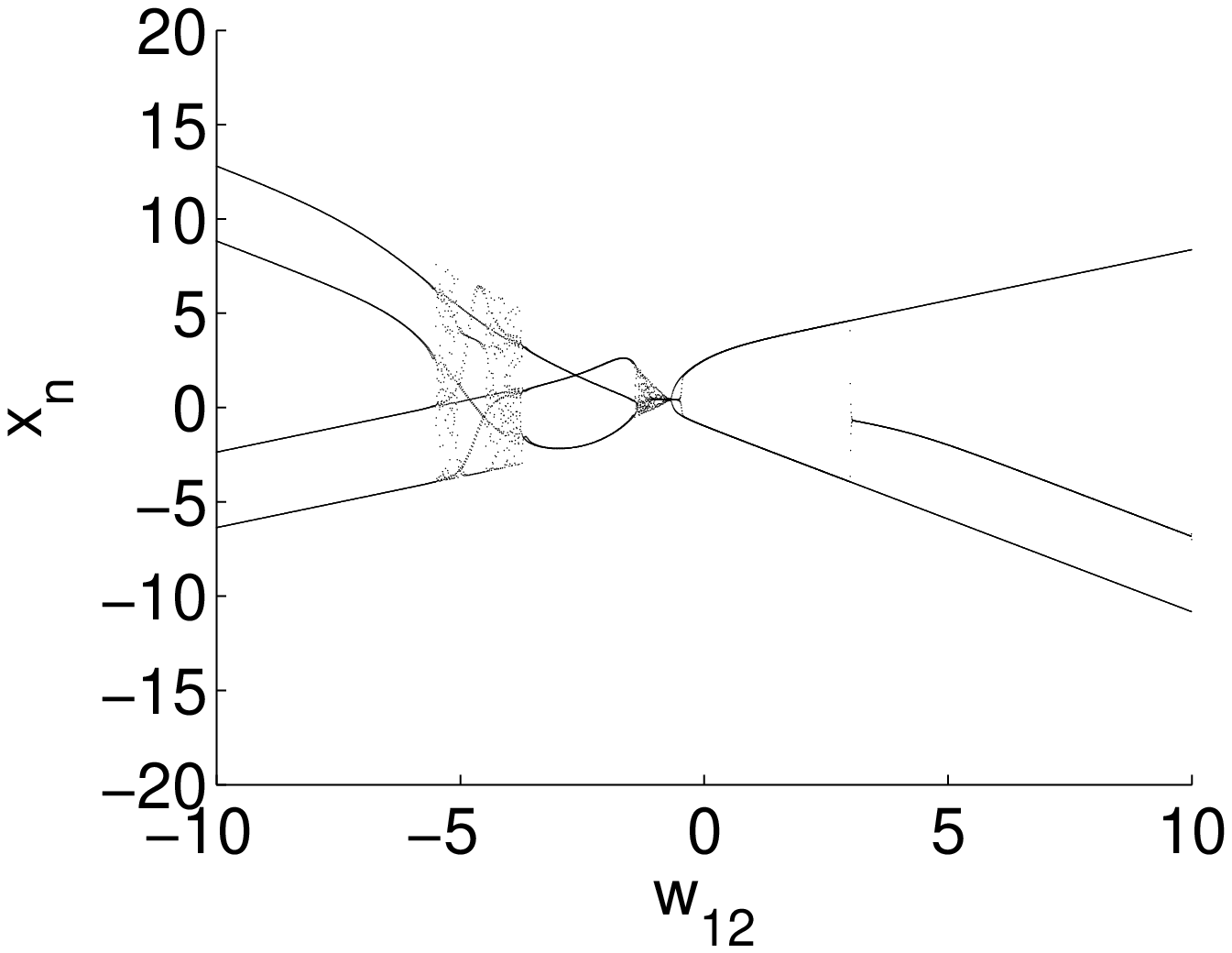} \hspace{0.2cm} 
\includegraphics[width=75mm]{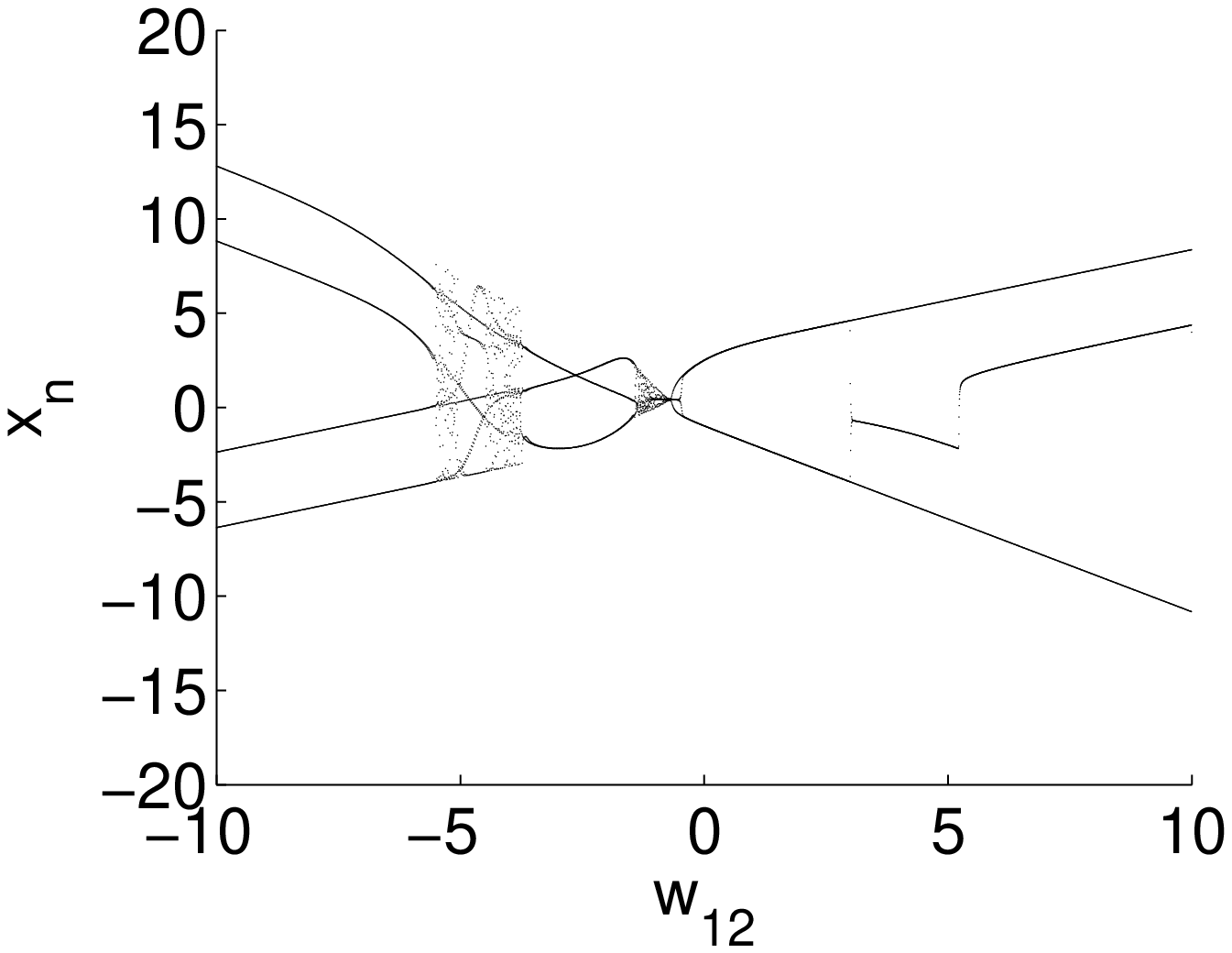}} 
\begin{tabbing}
\hspace*{3cm} \= \hspace*{8cm} \= \hspace*{2cm} \kill
 \> (a) \> (b) 
\end{tabbing}
\caption{\small{Typical bifurcation diagrams for system (\ref{eq:1}) when $w_{11}<0$. In this case, $w_{11}=-2$, %%@
$b_2=-3$, $w_{21}=5$, $\alpha=1$ and $\beta=0.3$. The parameter $w_{12}$ is ramped down and up, $-10 \leq w_{12} \leq %%@
10$, $b_1=1$. (a) The initial conditions are $x(0)=-7$, $y(0)=-7$. (b)  The initial conditions are $x(0)=4$, %%@
$y(0)=2$.}}
\end{figure}

\begin{figure}[htbp]
\centerline{\includegraphics[width=80mm]{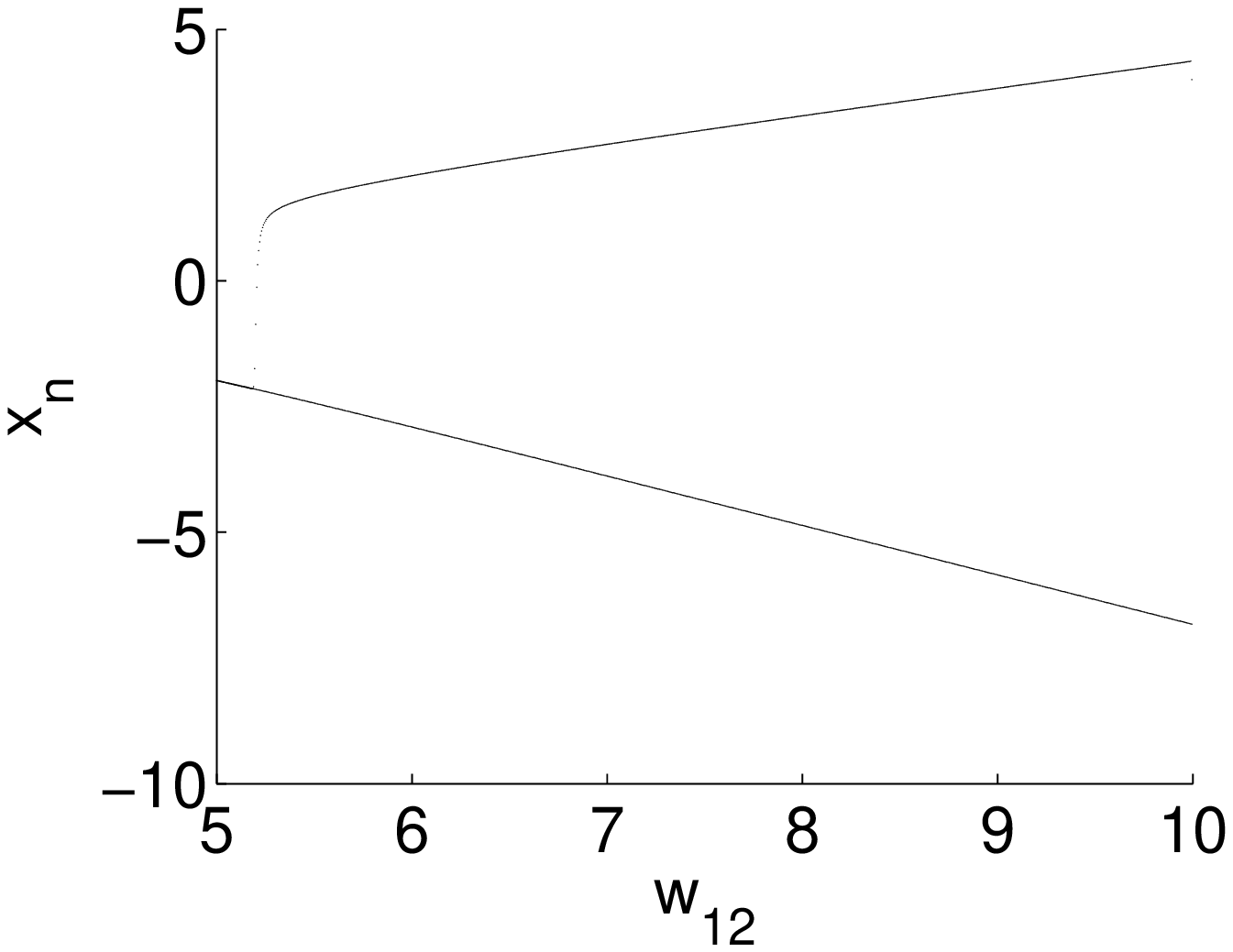}} 
\caption{\small{A bifurcation diagram for system (\ref{eq:1}) for the parameter values $w_{11}=-2$, $b_2=-3$, %%@
$w_{21}=5$, $\alpha=1$ and $\beta=0.3$. In this case, $w_{12}$ is ramped down from $w_{12}=10$ to $w_{12}=5$ and then %%@
ramped back up again. The initial conditions are $x(0)=4$, $y(0)=2$. There is a large open bistable region for %%@
$w_{12}>5$, approximately.}}
\end{figure}

In the final case, consider system (\ref{eq:2}) with the parameters listed as in Figure 10 and ramp $w_{12}$ down then up for $5 \leq w_{12} \leq 10$. The initial conditions chosen are $x(0)=4$ and $y(0)=2$ in order that a hysteresis cycle will exist, as Figure 13 demonstrates. There is a large open counterclockwise hysteresis cycle. If initial conditions $x(0)=-7$ and $y(0)=-7$ are chosen, then no bistable cycle will form. Thus, once more the system is shown to be history dependent.  \\

\noindent
{\bf 4. Conclusions and further work}
\vspace{0.2cm}

A stability analysis and feedback mechanism has been implemented to model the dynamics of neurons and simple %%@
neuromodules. Feedback is important in artificial neural networks, and is undoubtedly involved in the development of %%@
human brains. The dynamics of neurons and simple neuromodules demonstrate a broad range of possible behaviours for %%@
these systems. It has been shown that the dynamics are dependent not only on the history, but also on the initial %%@
conditions chosen at the start of the simulations, the step-size used in the simulations, the number of varying %%@
parameters and the extent to which those parameters are varied. This work clearly demonstrates that the stability %%@
analysis and the simultaneous plotting of bifurcation diagrams with feedback are required to fully understand the %%@
dynamics involved in these simple systems.

It is well known that stochastic resonance can cause jumps between attractors in the human brain. There appears to be %%@
little research on how history can affect jumps between different steady states in the brain. This work indicates the %%@
need to investigate feedback and its affect on neuronal models.

Future work will concentrate on feedback in highly nonlinear systems, where there are expected to be bifurcations %%@
involving multiple---fixed points, limit cycles and chaotic attractors. \\

\end{document}